\documentclass[12pt]{article}

\usepackage{a4wide, amsmath,amssymb,amscd,amsthm,makeidx,txfonts,graphicx,url,bm}

\usepackage[osf]{Baskervaldx} 
\usepackage{tkz-graph, pgfplots}

\usepackage{authblk}

\usepackage{color}
\definecolor{orange}{rgb}{.9,.6,.2}
\definecolor{violet}{rgb}{.5,0,.5}
\definecolor{dg}{rgb}{0,0.67,0}

\definecolor{cof}{RGB}{219,144,71}
\definecolor{pur}{RGB}{186,146,162}
\definecolor{greeo}{RGB}{91,173,69}
\definecolor{greet}{RGB}{52,111,72}

\renewcommand{\bar}{\overline}

\newcommand{\nc}{\newcommand}

\nc{\noi}{\noindent}
\nc{\cmmt}[1]{}

\newcounter{projectcnt}
\nc{\sctn}[1]{{\bigskip\addtocounter{projectcnt}{1}%
\begin{center}{\textbf{Project \theprojectcnt. #1}}\end{center}\medskip}}
\nc{\plainsctn}[1]{{\bigskip\begin{center}{\textbf{#1}}\end{center}\medskip}}

\nc{\dave}[1]{\begin{quote}\em #1\end{quote}}
\nc{\bbP}{{\bf{P}}}
\nc{\N}{{\bf{N}}}
\nc{\bF}{{\bar{\bf{F}}}}
\nc{\Gdave}{{\bf{G}}}
\nc{\I}{{\bf{I}}}
\nc{\cM}{{\cal M}}
\nc{\cK}{{\mathcal K}}
\nc{\cA}{{\cal A}}
\nc{\cO}{{\mathcal O}}
\nc{\cS}{{\mathcal S}}
\nc{\cNF}{{\cal NF}}
\nc{\cMF}{{\cal MF}}
\nc{\cLS}{{\cal LS}}
\nc{\Aff}{\mbox{Aff}}
\nc{\Mod}{{\mbox{\textup{Mod}}}}
\nc{\T}{\mathbb{T}}
\nc{\nf}{S_k^{\textup{\tiny new}}(N, \chi)}
\nc{\of}{S_k^{\textup{\tiny old}}(N, \chi)}
\nc{\mf}{S_k(N, \chi)}
\nc{\Skip}{\mathrm{\rm Skip}}
\nc{\Eis}{\mathrm{Eis}}
\nc{\Poly}{\mathrm{\rm Poly}}
\nc{\Polys}{\mathrm{\rm Polys}}
\nc{\Ext}{\mathrm{\rm Ext}}
\nc{\Ore}{\mathrm{\rm Ore}}
\nc{\Step}{\mathrm{\rm Step}}
\nc{\Mass}{\mathrm{\rm Mass}}
\nc{\mass}{\mathrm{\rm mass}}
\nc{\Vol}{\mathrm{\rm Vol}}

\newcommand{\m}{{\mathfrak m}}

\def\sump{\operatornamewithlimits{\sum\raise7pt\hbox{$\prime$}}}

\DeclareMathOperator{\id}{id}

\DeclareMathOperator{\charpoly}{char}
\DeclareMathOperator{\vol}{vol}

\DeclareMathOperator{\sign}{sign}

\def\sumpp{\operatornamewithlimits{\sum\raise9pt\hbox{\kern-2pt\scriptsize$(p)$\kern-11pt}\kern6pt}}

\DeclareMathOperator{\GL}{GL}

\DeclareMathOperator{\codeg}{codeg}

\newcommand{\calH}{{\cal H}}

\newcommand{\Q}{{\mathbb Q}}

\newcommand{\Z}{{\mathbb Z}}

\newcommand{\F}{{\mathbb F}}

\newcommand{\C}{{\mathbb C}}
\renewcommand{\P}{{\mathbb P}}

\newcommand{\A}{\mathbb A}

\newtheorem{theorem}{Theorem}[section]

\newcommand{\litem}{\par\noindent\dimen0=\parindent%
    \advance\dimen0 by-4pt
               \hangindent=\dimen0\ltextindent}

\newcommand{\ltextindent}[1]{\hbox to \hangindent{#1\hss}\ignorespaces}
\newcommand{\ltextjndent}[1]{\hbox to \hangindent{#1\hss}\ignorespaces\kern-1ex}

\begin{document}

\title{Mixed Hodge numbers and factorial ratios}

\author{Fernando Rodriguez Villegas 
\thanks{
I would like to thank A.~Mellit for several useful discussions on the
topic of this note}}
\affil{Abdus Salam Center for Theoretical Physics\\
{\tt villegas@ictp.it}}

\maketitle

\begin{abstract}
  This note is an extended version of the slides for my talk with the
  same title at the {\it Arithmetic, geometry, and modular forms: a
    conference in honour of Bill Duke} in June 2019 at
  the ETH in Z\"urich. The results presented concern three geometric
  criteria for the integrality of factorial ratios, numbers such as
  $(30n)!n!/(6n)!(10n)!(15n)!$,
  which are integral in a non-immediate way for all $n$.
  This work is an offshoot of an ongoing project on hypergeometric
  motives joint with D.~Roberts and M.~Watkins.
\end{abstract}

\section{}
\label{chebyshev}
In his work on the prime number theorem Chebyshev~\cite{Chebyshev}
used indirectly that the numbers
$$
c_n:=\frac{(30n)!n!}{(6n)!(10n)!(15n)!},
$$
are integral for $n=0,1,2,\ldots$. The only know proof of this fact
seems to be to use that the valuation
$$
v_p(c_n)\geq 0, \qquad n \geq 0
$$
for all primes $p$; this in turn relies on the fact that
$$
\lfloor 30x\rfloor +\lfloor x\rfloor -\lfloor 6x\rfloor -\lfloor 10x
\rfloor -\lfloor 15x\rfloor \geq 0
$$
for all real numbers $x$, which is easy to check. We will encode the
data defining $c_n$ by the list $\gamma:=(-30,-1,6,10,15)$ of positive
integers with zero sum.

I have been fascinated by this fact since pointed out to me by
P. Sarnak in the late 90's. It is not entirely clear how to understand
the integrality of ratios of factorials such as these
(see~\cite{Bober},~\cite{Borisov},~\cite{Sound},~\cite{Adolphson-Sperber}
and the bibliography therein for recent developments). In this note I
will discuss three different criteria. These involve

\bigskip
{\bf (A)} Interior lattice points of dilations of an associated
polytope~$\Delta$

\bigskip
{\bf (B)} Hodge numbers of a certain hypersurface $Z_t\subseteq
\T$ of a torus  $\T:=(\C^\times)^d$

\bigskip
{\bf (C)}
Effective weight of the corresponding hypergeometric motive $\calH(\gamma\,|\,t)$

\bigskip
\section{} 
Before discussing these criteria, I would like to briefly sketch the
connection between the integrality of $c_n$ and the algebraicity of
the corresponding power series.

Early on I noticed that the hypergeometric series
\begin{equation}
\label{hyperg-series}
c(t):=\sum_{n\geq 0} c_n \left(\frac tM\right)^n =
{}_8F_7\left[
 \begin{matrix}
 \tfrac1{30}&\tfrac7{30}&\tfrac{11}{30}&\tfrac{13}{30}&\tfrac{17}{30}&
 \tfrac{19}{30}&\tfrac{23}{30}  &\tfrac{29}{30} \\
 \tfrac 12 & \tfrac13&\tfrac23&\tfrac15&\tfrac25&\tfrac35&\tfrac45&
\end{matrix}  \left.\right|\,t\right], \qquad
\qquad M:=\frac{30^{30}}{6^6{10}^{10}15^{15}},
\end{equation}
is an algebraic function of $t$.
This follows from the work of Beukers and Heckman~\cite{BH} who
classified all algebraic hypergeometric series, extending the
classical result of Schwartz on the algebraicity of the ${}_2F_1$
hypergeometric series. Our series $c$ is number $67$ in the
Beukers-Heckman list.

Concretely, there exists a polynomial $A(x,y)\in \Z[x,y]$, which we
may assume is irreducible, such that
$$
A(t,c(t))=0
$$
in the ring of power series in $t$. What is the connection between
the algebraicity of $c$ and the integrality of the coefficients $c_n$?

A little work using a theorem of Eisenstein shows that the latter is a
consequence of the former (see~\cite{FRV-fields}[Prop.2] for the
details). But algebraicity is a much stronger property of a power
series than integrality of its coefficients. Nevertheless, there is
some connection and Landau~\cite{Landau} had already exploited it to
reprove Schwartz's result.

The Galois group $\Gamma$
of the normal closure of the extension $\Q(c(t))/\Q(t)$ is
the Weyl group $W(E_8)$
of the Lie algebra $E_8$
of order $696,729,600$. It turns out, as we will see shortly, that
$$
\deg_y A= 483,840.
$$
So it is not going to be very easy to find $A$ (though
D.~Roberts~\cite{Roberts} has computed a degree $240$ polynomial with
the same Galois closure.) 

On the other hand, the series $c$
satisfies a linear differential equation with polynomial coefficients
of order $8$.
This equation has only regular singularities at $t=0,1,\infty$.
The corresponding {\it monodromy representation}, obtained by analytic
continuation of local solutions at some point $t_0\neq 0,1,\infty$,
$$
\rho: \quad
\pi_1\left(\P^1\setminus\{0,1,\infty\}\right)\longrightarrow \GL(V), 
$$
has the following properties:
$$
q_\infty:=\charpoly(T_\infty)= \Phi_{30},
\qquad
q_0:=\charpoly(T_0^{-1})= \Phi_{1}\Phi_{2}\Phi_{3}\Phi_{5}
$$
and $T_1$ fixes a codimension one subspace of $V$. Here $\Phi_n$
denotes the $n$th cyclotomic polynomial, $V$ is the
space of local solutions to the differential equation around $t_0$ and
$T_s$, for $s=0,1,\infty$, are the local monodromies; i.e., the images by
$\rho$ of small loops around $s$. These loops are chosen oriented in a
consistent manner so  that $T_0T_1T_\infty= \id_V$. The image of
$\rho$ is called the {\it monodromy group}. 

  Note that
$$
\frac{\Phi_{30}}{\Phi_{1}\Phi_{2}\Phi_{3}\Phi_{5}}=
\frac{(T^{30}-1)(T-1)}{(T^6-1)(T^{10}-1)(T^{15}-1)}
$$
connecting the list $\gamma$ determining $c$ and the hypergeometric
parameters in the series~\eqref{hyperg-series}. 

In general, the monodromy representation of a hypergeometric
differential equations is uniquely determined by the analogue of the
polynomials $q_\infty,q_0$ (when irreducible or equivalently when
$q_0$ and $q_\infty$ have no common roots). In other words, they give
rise to a rigid local system (in the sense of N. Katz~\cite{Katz1}):
the local monodromies uniquely determine the monodromy
representation. Clasically, this is known as not having {\it accessory
  parameters}.

We can build this representation in our case by starting directly with
the group $W:=W(E_8)$
(see~\cite{Bourbaki} for details on Coxeter groups).  A Coxeter
element $\sigma\in W$
is the product of all simple reflections. Its conjugacy class is
uniquely determined independent of the order in which we perform the
product. The order of $\sigma$
is the Coxeter number, which equals $30$ for $E_8$.

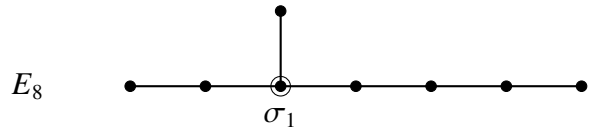
\begin{center}
\begin{tikzpicture}
    \draw (-1,0) node[anchor=east]{$E_8$};

    \draw[fill=black] ( 0, 0) circle (.07);
    \draw[fill=black] ( 1, 0) circle (.07);
    \draw[fill=black] ( 2, 0) circle (.07);
    \draw[fill=black] ( 3, 0) circle (.07);
    \draw[fill=black] ( 4, 0) circle (.07);
    \draw[fill=black] ( 5, 0) circle (.07);
    \draw[fill=black] ( 6, 0) circle (.07);
    \draw[fill=black] ( 2, 1) circle (.07);

    \draw[thick] (0,0) -- (6,0);
    \draw[thick] (2,0) -- (2,1);

    \draw (2,0) circle (.13);
    \node  [label=below:{$\sigma_1$}] at (2,0) {};
\end{tikzpicture}
\end{center}

Take $\tau=\sigma_1\cdots\sigma_8$. Then
$\sigma_1\tau=\sigma_2\cdots\sigma_8$ is the product of the Coxeter
elements of the diagram obtained by removing the circled dot and all
its attached edges in the Dynkin diagram of $E_8$. Its characteristic
polynomial is then $q_0$.

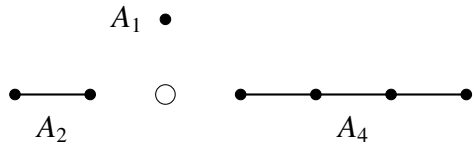
\begin{center}
\begin{tikzpicture}
    \draw (-1,0) node[anchor=east]{};

    \draw[fill=black] ( 0, 0) circle (.07);
    \draw[fill=black] ( 1, 0) circle (.07);
    \draw[fill=black] ( 3, 0) circle (.07);
    \draw[fill=black] ( 4, 0) circle (.07);
    \draw[fill=black] ( 5, 0) circle (.07);
    \draw[fill=black] ( 6, 0) circle (.07);
    \draw[fill=black] ( 2, 1) circle (.07);

    \draw[thick] (0,0) -- (1,0);
    \draw[thick] (3,0) -- (6,0);

   \node  [label=left:{$A_1$}] at (2,1) {};
   \node  [label=below:{$A_2$}] at (0.5,0) {};
   \node  [label=below:{$A_4$}] at (4.5,0) {};

    \draw (2,0) circle (.13);
\end{tikzpicture}
\end{center}

It follows that if we choose
$$
T_\infty := \tau, \qquad T_1:=\sigma_1, \qquad
T_0^{-1}:=\sigma_1\tau
$$
we obtain a representation isomorphic to $\rho$ by rigidity. It is
straightforward to check using a computer that the group generated by
$T_0,T_1,T_\infty$, the monodromy group, is all of $W(E_8)$. 

In our case the Galois group $\Gamma$
coincides with the monodromy group. The degree of $A$
in $y$
can now be computed by Galois theory. Indeed, using basic properties
of reflection groups we find that it equals
$$
483,840=[W(E_8):W(A_1)\times W(A_2)\times W(A_4)]
$$
as promised. (It suffices to find the simple roots orthogonal to that
corresponding to $\sigma_1$.)

\section{}
We now turn to criterion~(A) and the associated polytope.  Let
$\gamma=(\gamma_1,\ldots,\gamma_l)$ be a non-empty list of non-zero
integers with zero sum, no pair of entries satisfying
$\gamma_i+\gamma_j=0$ and with no positive integer dividing every
entry. Let $r$ be the number of negative $\gamma_i$ and $s$ the number
of positive ones so that $r+s=l$. Assume further that $s\geq r$. We
will call such a list a {\it gamma list} for short. Since the order of
the $\gamma_i$'s is irrelevant we will typically choose
$\gamma_1\leq \gamma_2\leq \cdots \leq \gamma_l$.  The numbers we
would like to study are then
$$
c_n:=\frac{\prod_{\gamma_i <0}(-\gamma_in)!}{\prod_{\gamma_i >
    0}(\gamma_i  n)!}, \qquad n=0,1,2,\ldots
$$

Given a gamma list consider $m_1,\ldots, m_l\in
\Z^d$ with $d:=l-2$ such that $\gamma$ spans their affine
relations:
$$
\gamma_1m_1+\cdots  \gamma_l m_l =  0.
$$
We will choose the $m_i$'s so that their affine span is primitive. In
practice, to find $m_i$ we can simply drop say $\gamma_1$ from the
list and find generators $m_2,\ldots,m_l$ over $\Z$ of the kernel of
the $1\times d$ matrix $(\gamma_2,\ldots,\gamma_l)$ and then set
$m_1=0$. Finally, we let $\Delta\subseteq \Z^d$ be the convex hull of
$m_1,\ldots, m_l$.  Primitivity guarantees that $\Delta$ is uniquely
determined up to invertible affine linear transformation over~$\Z$.
 The normalized volume of $\Delta$ equals 
$$
\vol(\Delta)=\vol(\gamma):=-\sum_{\gamma_i<0}\gamma_i=\sum_{\gamma_i>0}\gamma_i. 
$$

A possible choice of $m_i$ in the Chebyshev case are the columns of
the following matrix
$$
\left(
\begin{array}{ccccc}
1&0&5&0&0\\
1&0&0&3&0\\
1&0&0&0&2
\end{array}
\right)
$$
Here is a schematic (and not to scale) picture of the polytope with
each $m_i$ labeled by the corresponding $\gamma_i$.

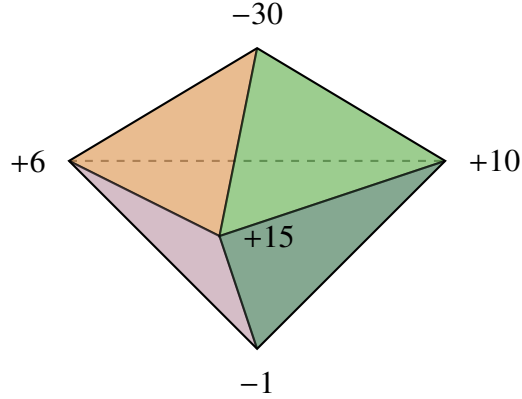
\begin{center}
\begin{tikzpicture}[thick,scale=5]
\coordinate (A1) at (0,0);
\coordinate (A3) at (1,0);
\coordinate (A4) at (0.4,-0.2);
\coordinate (B1) at (0.5,0.3);
\coordinate (B2) at (0.5,-0.5);

\begin{scope}[thick,dashed,,opacity=0.6]
\draw (A1) -- (A3);
\end{scope}
\draw[fill=cof,opacity=0.6] (A1) -- (A4) -- (B1);
\draw[fill=pur,opacity=0.6] (A1) -- (A4) -- (B2);
\draw[fill=greeo,opacity=0.6] (A3) -- (A4) -- (B1);
\draw[fill=greet,opacity=0.6] (A3) -- (A4) -- (B2);
\draw (B1) -- (A1) -- (B2) -- (A3) --cycle;

\node  [label=left:{$+6$}] at (A1) {};
\node  [label=right:{$+10$}] at (A3) {};
\node  [label=right:{$+15$}] at (A4) {};
\node  [label=above:{$-30$}] at (B1) {};
\node  [label=below:{$-1$}] at (B2) {};
\end{tikzpicture}
\end{center}

In general, for a $d$-dimensional  polytope $\Delta\in \Z^d$  we have, by
results of  Ehrhart,
\begin{equation}
\label{ehrhart}
\sum_{k\geq 0} \#(k\Delta)\, T^k =  \frac{\delta(\Delta,T)}{(1-T)^{d+1}},
\end{equation}
with $\delta(\Delta,T)$ a polynomial of degree at most $d$. We define
the {\it codegree} of $\Delta$ as
$\codeg(\Delta):=d+1-\deg(\delta)\geq 1$. Equivalently, the codegree
is the smallest positive integer $k$ such that $k\Delta$ has an
interior lattice point (see~\cite{Beck} for a recent survey of Ehrhart
theory). 

So for the Chebyshev example we have
$$
\sum_{k\geq 0} \#(k\Delta)\, T^k =  \frac{1+15T+15T^2}{(1-T)^4} =
1+19T+85T^2+230T^3+O(T^4)
$$
and hence $\codeg(\Delta)=2$.

We can now formulate our first criterion for integrality (see
also~\cite{Adolphson-Sperber}). 

\bigskip
\noindent
{\bf Criterion A}:
\begin{theorem}
\label{integr-1}
The numbers $c_n\in \Z$ for all $n$ if and only if  $s> r$  and 
\begin{equation}
\label{codegree}
\codeg(\Delta)\geq r;
\end{equation}
equivalently, $k\Delta$ has no interior lattice points for
$k=1,2,\ldots, r-1$.
\end{theorem}

Note that for $r=1$ the condition~\eqref{codegree} is vacuous. Indeed, 
in this case $c_n$ is a multinomial number and integrality is
immediate.

\section{}
We now define the algebraic varieties $Z_t$
appearing in criterion~(B) (see the MAGMA manual's {\it canonical
  scheme}~\cite{MAGMA} and~\cite{BCM}).  Consider the Laurent
polynomial, with $m_1,\ldots,m_l$ as before,
$$
f:=\sum_{i=1}^l u_ix^{m_i},\qquad x=(x_1,\ldots,x_d),
$$
where $x^m:=x_1^{m^1}\cdots x_l^{m^l}$ for
$m=(m^1,\ldots,m^l)\in \Z^d$ and $u_1,\ldots, u_l$ are
parameters in $\C^\times$. Following Watkins and Beukers-Cohen-Mellit
we consider the hypersurface $Z\subseteq \T:=(\C^\times)^d$ for a
given $u$ defined by the vanishing of~$f$. The dimension of $Z$ is
$\kappa:=l-3$.

By scaling the variables by $x_j\mapsto a_jx_j$ and the polynomial
itself by $f\mapsto a_0f$ with $a_0,a_j \in \C^\times$ for
$j=1,\ldots, d$ we obtain isomorphic hypersurfaces. Hence we may take
as the natural parameter for the family the quantity
$u:=u_1^{\gamma_1}\cdots u_l^{\gamma_l}\in \C^\times$. We see that
having choosen $f$ to have $l=d+2$ monomials in $d$ variables is what
guarantees that our family really only depends on one parameter.

In fact, it is better to normalize the parameter and use instead 
$$
t:=(-1)^{\vol(\gamma)}uM, \qquad M:= \frac{\prod_{\gamma_i
    <0}(-\gamma_i)^{-\gamma_i}}{\prod_{\gamma_i > 
    0}\gamma_i^{\gamma_i} }.
$$

Choose a family of hypersurfaces $Z_t$ with given parameter
$t\in \C^\times$. Concretely, we may choose integers $k_1,\ldots, k_l$
such that $k_1\gamma_1+\cdots k_l\gamma_l=1$ and take $u_i=u^{k_i}$
with $u=(-1)^{\vol(\gamma)}t/M$.  Then $Z_t$ is smooth except for
$t=1$ where it has a unique double point (uniqueness follows from the
primitivity of $m_1,\ldots,m_l$).

There is a refinement $\delta^\#$
of the polynomial $\delta$
in~\eqref{ehrhart} due to Stanley~\cite{Stanley} that incorporates the
way the faces of $\Delta$
are put together. By work of Danilov-Khovanski, Batyrev-Borisov,
Katz-Stapledon and others (see~\cite{DK}, \cite{BB}, \cite{KS}) this
polynomial gives the weight $\kappa$
mixed Hodge numbers~\cite{Deligne-1},~\cite{Deligne-2} of the middle
cohomology of $Z_t$ for generic $t\in \C^\times$. More precisely,
$$
\delta^\#(\Delta,T)=\sum_{j=0}^\kappa h_c^{j,\kappa-j,\kappa}(Z_t)T^{j+1}. 
$$

Since $\Delta$ is a simplicial polytope (all proper faces are
simplices) it is fairly straightforward to compute this polynomial
explicitly and obtain 
\begin{equation}
\label{stanley}
\delta^\#(\Delta,T)=\sum_{\substack{N\geq 1\\m_->m_+}}
\frac{T^{m_-}-T^{m_+}}{T-1}\delta_N^\#(T), \qquad \qquad
\delta_N^\#(T):=\sum_{\substack{j=1\\\gcd(j,N)=1}}^{N-1}
T^{\sum_{i=1}^l\left\{\frac{j\gamma_i}N\right\}}.
\end{equation}
Here 
$$
m_{\pm}:=\{i \,|\,  \sign(\gamma_i)= \pm1, N\mid \gamma_i\}
$$
and $\{\cdot\}$ denotes fractional part. One may verify that this
formula is equivalent to that giving the Hodge numbers of
hypergeometric motives first conjectured by Corti and Golyshev and
proved by Fedorov (see~\cite{CG}, \cite{HGM-Book}
and~\cite{Fedorov}). A completely analogous formula holds where we sum
over all $N$ such that $m_- < m_+$. In fact there is also a similar
formula for the polynomial $\delta(\Delta,T)$ itself.
$$
\delta(\Delta,T)=\sum_{N\geq 1}
\frac{T^{m_-}-1}{T-1}\delta_N^\#(T)
=\sum_{N\geq 1}
\frac{T^{m_+}-1}{T-1}\delta_N^\#(T).
$$

For our running example of Chebyshev with $\gamma=(-30,-1,6,10,15)$, 
we find that the only contribution to the sum comes from $N=30$
for which $\delta_{30}^\#(T)=8T^2$. Hence also
$\delta^\#(\Delta,T)=8T^2$. We can take
$$
Z_t: \quad xyz  -\frac Mt + x^5 + y^3 + z^2=0,\qquad \qquad
M:=\frac{30^{30}}{6^610^{10}15^{15}},
$$
an affine piece of a rational elliptic surface with parameter $x$; it
has twelve bad fibers of type $I_1$.  We have $h_c^{i,j,2}(Z_t)=0$ for
$(i,j)=(2,0),(0,2)$ and $h_c^{i,j,2}(Z_t)=8$ for $(i,j)=(1,1)$. This
exhibits the connection with the $E_8$ lattice more directly as the
Mordell-Weil lattice of the elliptic surface~\cite{schuett-shioda}.

In general, $\deg \delta^\#\leq \deg \delta$ but for our polytopes the
equality holds. We may therefore reformulate Theorem~\ref{integr-1} as
follows.

\bigskip
\noindent
{\bf Criterion B}:
\begin{theorem}
\label{integr-2}
With the notation of Theorem~\ref{integr-1} we have 
that $c_n\in \Z$ for all $n$ if and only if  $s> r$  and 
\begin{equation}
\label{hodge}
h^{\kappa,0}(Z_t)=h^{\kappa-1,1}(Z_t)=\cdots=h^{\kappa-r+1,r-1}(Z_t)=0,
\qquad t \neq 1.
\end{equation}
\end{theorem}
For example, for $r=2$ the condition~\eqref{hodge} means that
the hypersurface $Z_t$ should have geometric genus $p_g:=h^{\kappa,0}$
equal zero.

In fact, we can go further and compute all of the mixed Hodge numbers
of $Z_t$
using a formula of Batyrev-Borisov~\cite{BB}[Thm. 3.18] and not just
the top degree piece given by~\eqref{stanley}. It is better to
formulate the result in terms of the {\it primitive} cohomology
$PH_c^\kappa(Z_t)$;
i.e., the kernel of the Gysin homomorphism~\cite{DK}[Prop. 3.9]
$H_c^i(Z_t,\C)\rightarrow H_c^{i+2}(\T,\C)$.
This removes from the cohomology of $Z_t$
the contributions from the ambient torus.

\begin{theorem} The $E$-polynomial of
$PH_c^\kappa(Z_t)$ (primitive cohomology) is
$$
E(\Delta;a,b):=\sum_{i,j}h_c^{i,j}a^ib^j=
\frac1{ab}\left[\delta^\#(\Delta;a,b)+
\delta^0(\Delta;a,b)-1\right],
$$
where
$$
\delta^\#(\Delta;a,b)=\sum_{\substack{N\geq 1\\m_->m_+}}
\frac{(a/b)^{m_-}-(a/b)^{m_+}}{a/b-1} b^{l-1}\delta^\#_N(a/b)
$$
$$
\delta^0(\Delta;a,b)=\sum_{N\geq 1}
\frac{(ab)^{\min(m_-,m_+)}-1}{ab-1} b^{l-m_+-m_-}\delta^\#_N(a/b)
$$
\end{theorem}

Note that $T\delta^\#(\Delta,T)=T\delta^\#(\Delta;T,1)$ is indeed the
component of (top) degree $\kappa$ of $E(\Delta;a,b)$. Also
$$
\delta(\Delta,T)=TE(\Delta;T,1)+1.
$$
It is useful to display the two variable polynomial $E(\Delta;a,b)$ as
an array with the coefficient of $a^ib^j$ in spot $(i,j)\in \Z^2$. For
Chebyshev we have
$$
\begin{array}{c|c|c|l}
N& m_-&m_+&\delta^\#\\
\hline
1& 2& 3& 1\\
2& 1& 2& T\\
3& 1& 2& 2T\\
5& 1& 2& 4T\\
6& 1& 1& T^2+T\\
10& 1& 1& 2T^2+2T\\
15& 1& 1& 4T^2 + 4T\\
30&  1& 0&  8T^2
\end{array}
$$
and hence
$$
E(\Delta;a,b)=
\begin{array}{cc}
7&8\\
8&7
\end{array}
$$

\section{}

There is also an arithmetic aspect to the story. The number of points
of $Z_t$
over finite fields (\cite{BCM}) involves the finite analogue of a
hypergeometric series due to N. Katz~\cite{Katz}[(8.2.7)]. More
precisely, the trace of the $q$-th
Frobenius (for a good prime $p$)
acting on the weight $\kappa$ piece of the middle cohomology of $Z_t$
has the form
\begin{equation}
\label{hyperg-trace}
\calH(t):=\frac 1{1-q}\sum_\chi
\frac{J(\alpha\chi,\beta\chi)}{J(\alpha,\beta)}\chi(t),
\end{equation}
where $\alpha,\beta$ are the hypergeometric parameters (as characters
of $\F_q^\times$), $J$ are certain Jacobi sums and $\chi$ runs over
all characters of $\F_q^\times$. We call the weight $\kappa$ piece of the
middle cohomology the {\it hypergeometric motive}
$\calH(\gamma|t)$~\cite{HGM-Book} attached to the list $\gamma$. It is
a pure motive of weight~$\kappa$.  

For our running Chebyshev example, we consider the hypersurface of
$\A^3$ given by
$$
\tilde Z_t: \quad -\tfrac Mt+xyz+x^5+y^3+z^2=0, \qquad
M:=\frac{30^{30}}{6^610^{10}15^{15}},
$$
and find (see~\cite{BCM}[Cor. 1.8] for an equivalent statement) that 
$$
\#\tilde Z_t(\F_q) = q^2 +q\calH(t) +1, \qquad t\neq 0, \quad
q=p^k,\quad  p>5. 
$$
A number of features of hypergeometric motives have already been
implemented in MAGMA by M. Watkins~\cite{MAGMA}. For example, we can
easily verify a few instances of the above formula as follows.

Define the function
\begin{verbatim}
function surfcheb(t,q);

         A<x,y,z>:=AffineSpace(GF(q),3);
         M:=30^30/6^6/10^10/15^15;u:=-M/t;
         S:=Surface(A,u+x*y*z+z^5+x^3+y^2);
         return((#Points(S)-q^2-1)/q);

end function;
\end{verbatim}
This will output the number $(\#\tilde Z_t(\F_q)-q^2-1)/q$ for a given
value of $t$ and $q$. Then for example 
\begin{verbatim}
> H:=HypergeometricData([*-30,-1,6,10,15*]);
> [HypergeometricTrace(H, t, 7): t in [1..6]];
[ 0, 0, 1, -1, -1, 0 ]
> [surfcheb(t,7): t in [1..6]];             
[ 0, 0, 1, -1, -1, 0 ]
> [HypergeometricTrace(H, t, 23): t in [1..22]];
[ 0, 0, 1, -1, 0, 0, -1, 2, -1, 0, 1, 0, -2, -2, 0, 1, 0, 0, 1, 0, -1, 1 ]
> [surfcheb(t,23): t in [1..22]];            
[ 0, 0, 1, -1, 0, 0, -1, 2, -1, 0, 1, 0, -2, -2, 0, 1, 0, 0, 1, 0, -1, 1 ]
\end{verbatim}
In fact, the MAGMA package can compute many more things about
$\calH(\gamma|t)$, notably, a big chunk (and in many cases all) of its
$L$-function. I refer to the MAGMA manual~\cite{MAGMA} for many worked out
examples. Here are for example, some Euler factors for
$\calH(\gamma\,|\,2)$. 
\begin{verbatim}
> EulerFactor(H,2,7);
x^8 + x^5 + x^3 + 1
> EulerFactor(H,2,23);
x^8 - x^4 + 1
\end{verbatim}

We can formulate our last criterion in the following form.

\bigskip
\noindent
{\bf Criterion C}:
\begin{theorem}
\label{hgm}
  With the notation of Theorem~\ref{integr-1} we have that $c_n\in \Z$
  for all $n$ if and only if $s> r$ and the hypergeometric motive
  $\calH(\gamma|t)$ is a Tate twist of a pure effective motive of
  weight   $s-r-1$. 
\end{theorem}
For example, if $s-r=1$
then the integrality of $c_n$
for all $n$
is equivalent to $\calH(\gamma|t)$
being a Tate twist of a pure motive of weight zero. It can be shown
directly~\cite{FRV-fields}[Thm. 1.4] that when $s-r=1$
the integrality of $c_n$
is equivalent to Beukers and Heckman's criterion for the algebraicity
of the power series $c=\sum_{n\geq 0} c_n t^n$.
Hence in a sense the above theorem is a generalization of the this
criterion.

For $s-r=2$ the integrality of $c_n$ is equivalent to the hypergeometric
motive $\calH(\gamma|t)$ being a Tate twist of a motive with Hodge
numbers $(m,m)$ for some positive integer $m$. Conjecturally,
$\calH(\gamma|t)$ should then correspond to an abelian variety.

For example, consider $\gamma=(-11,-2,1,3,4,5)$.
The characteristic polynomials of local monodromies are
$$
q_\infty=\Phi_{11},\qquad \qquad 
q_0=\Phi_1\Phi_3\Phi_4\Phi_5.
$$
We can take $\Delta$ as the convex hull of
the columns of the matrix
$$
\left(
\begin{array}{cccccc}
1&0&0&1&2&0\\
1&0&0&2&0&1\\
1&0&1&0&0&2\\
1&1&2&0&0&0
\end{array}
\right)
$$
and $Z_t$ as the family of cubic threefolds with equation
$$
Z_t:\quad x_1^2+x_1x_2^2+x_2x_3^2-\tfrac tM x_3x_4^2+x_4^2+x_1x_2x_3=0, \qquad
t\in \C^\times,  \quad  M:=\frac{11^{11}2^2}{3^34^45^5}.
$$
Taking the Zariski closure $\overline Z_t \subseteq \P^4$ we obtain a
family of projective cubic threefolds, smooth for $t\neq 1$. It is a
classical fact that the Hodge numbers of a smooth projective cubic threefold
for its middle cohomology are
$$
(h^{3,0},h^{2,1},h^{1,2},h^{0,3})=(0,5,5,0).
$$ 
This matches the calculation of $\delta^\#(\Delta,T)$ using~\eqref{stanley}
and the values below.
$$
\begin{array}{c|c|c|l}
N& m_-&m_+&\delta_N^\#\\
\hline
1 & 2& 4&  1\\
2 & 1& 1&  T^2\\
3 & 0& 1&  T^3 + T^2\\
4 & 0& 1&  T^3 + T^2\\
5 & 0& 1&  2T^3 + 2T^2\\
11& 1& 0&  5T^3 + 5T^2\\
\end{array}.
$$
We also obtain $\delta(\Delta,T)=5T^3 + 6T^2 + T + 1$ and
$$
E(\Delta;a,b) =
\begin{array}{ccc}
0 &5&0\\
0& 1& 5\\
1& 0& 0
\end{array}
$$

Here the hypergeometric motive $\calH(\gamma | t)$ is a Tate twist of
$H^1$ of the intermediate Jacobian of $\overline Z_t$
(see~\cite{Clemens-Griffiths} and~\cite{Bombieri-Swinnerton-Dyer}) and
$$
\#\overline Z_t(\F_q) = q^3+q^2+q+1-q\calH(t), \qquad \qquad
q=p^k,\quad p >11.
$$
We can independently check that consistent with Theorem~\ref{hgm} the
following factorial ratios
$$
c_n:=\frac{(11n)!(2n)!}{n!(3n)!(4n)!(5n)!},
$$
are integral for $n=0,1,\ldots$

As a final example consider
$$
\gamma:=(-63,-8,-2,1,4,16,21,31).
$$
Now we have
$$
q_\infty= \Phi_{63}\Phi_3, \qquad \qquad
q_0=\Phi_1^2\Phi_4\Phi_{16}\Phi_{31}.
$$
We may take the projective closure of $Z_t$ in $\P^5$
$$
\overline Z_t:\quad x_1^2x_3+x_0x_1x_2+x_1x_2^2-\tfrac tM
x_2x_6^2+x_3x_5+x_4^2x_6+x_4x_5^2+x_0^3, \qquad t\in \C^\times,
$$
a smooth family of cubic fivefolds (for $t\neq 1$).
The only non-zero Hodge numbers of the middle cohomology of
$\overline Z_t$ are $h^{4,1}=h^{1,4}=21$ and
$\calH(\gamma|t)=H^5(\overline Z_t)$ is a pure motive of rank $42$ and
weight $5$.   

The data on the polytope $\Delta$ is
$$
\begin{array}{c|c|c|l}
N& m_-&m_+&\delta^\#\\
\hline
1& 3& 5& 1\\
2& 2& 2& T^2\\
3& 1& 1& T^4 + T^2\\
4& 1& 2& T^3 + T^2\\
7& 1& 1& 6T^3\\
8& 1& 1& 4T^3\\
9& 1& 0& 3T^4 + 3T^3\\
16& 0& 1& 4T^4 + 4T^3\\
21& 1& 1& 12T^3\\
31&  0& 1& 15T^4 + 15T^3\\
63&  1& 0&  18T^4 + 18T^3
\end{array}
$$
giving $\delta^\#(\Delta,T)=21T^4+21T^3$, $\delta(\Delta,T)=22T^4 +
45T^3 + 4T^2 + T + 1$ and
$$
E(\Delta;a,b)=
\begin{array}{cccc}
0& 1& 21&  0\\
0& 1& 23& 21\\
0& 2&  1&  1\\
1& 0&  0 & 0
\end{array}
$$

In this case,
$$
\#\overline Z_t(\F_q)=q^5+q^4+q^3+q^2+q+1-q^2\calH(t), \qquad
q=p^k,\quad p > 31.
$$
Again we find that
$$
c_n:=\frac{(63n)!(8n)!(2n)!}{n!(4n)!(16n)!(21n)!(31n)!},
$$
is integral for $n=0,1,\ldots$

\end{document}